\documentclass[12pt]{amsart}
\usepackage{amssymb,latexsym,amscd}

\numberwithin{equation}{section}

\newtheorem{theorem}{Theorem}[section]
\newtheorem{proposition}[theorem]{Proposition}

\newtheorem{lemma}[theorem]{Lemma}

\theoremstyle{definition}

\newtheorem{example}[theorem]{Example}
\newtheorem{definition}[theorem]{Definition}

\newtheorem{problem}[theorem]{Problem}

\renewcommand{\eqref}[1]{{\rm (\ref{#1})}}

\def\endproof{\hfill$\square$\medskip}

\def\AA{\mathbb{A}}

\def\ZZ{\mathbb{Z}}
\def\CC{\mathbb{C}}

\def\QQ{\mathbb{Q}}
\def\GG{\mathbb{G}}

\def\BB{\mathcal{B}}

\def\ii{\mathbf{i}}

\thanks{The authors were supported in part
by NSF grants  (A.B.), and  by ISF and NSF grants (D.K.).
}

\begin{document}

\title{Lecture notes on Geometric Crystals and their combinatorial analogues}

\author{Arkady Berenstein}
\address{Department of Mathematics, University of Oregon,
Eugene, OR 97403, USA} \email{arkadiy@math.uoregon.edu}

\author{David Kazhdan}
\address{\noindent Department of Mathematics, Heberew University, Jerusalem, Israel}
\email{kazhdan@math.huji.ac.il}

\date{November 21, 2005}

\makeatletter
\renewcommand{\@evenhead}{\thepage \hfill  ARKADY BERENSTEIN and DAVID KAZHDAN \hfill \thepage}

\renewcommand{\@oddhead}{\hfill Lecture notes on geometric crystals\hfill \thepage}
\makeatother

\maketitle

\tableofcontents

\section{First lecture: Geometric and unipotent crystals}
\medskip

The first author would like to express his gratitude to Professor Okado for the opportunity to give the lectures, to Professor Kuniba for the great help in preparation of the lecture notes, and to RIMS for the hospitality. 

The lectures are based on the results of our research originated six years ago in \cite{bk1} and continued in \cite{bk2}. We start with the problem of birational Weyl group actions that served as the original motivation for this work (for the terminology and results on reductive algebraic groups, see, e.g., \cite{ov}).

\begin{problem}  

\label{prob:birational Weyl group actions}

Let $G$ be a split reductive algebraic group with a maximal torus $T$. Given an affine algebraic variety $X$, a function $f$ on $X$,  and a morphism of algebraic varieties $\gamma:X\to T$, construct a birational action of the Weyl group $W=Norm_T(G)/T$ on $X$ in such a way that:

\smallskip 

\noindent (1) The structure map $\gamma:X\to T$ commutes with the $W$-action (where the $W$-action on $T$ is the natural one).

\noindent (2) The function  $f$ is $W$-invariant.

\noindent (3) For each $w\in W$ the fixed point set $X^w=\{x\in X: w(x)=x\}$ is the pre-image $\gamma^{-1}(T^w)$ of the fixed point set $T^w=\{t\in T: w(t)=t\}$ (i.e., all fixed points of $w$ ''upstairs'' come from the fixed points of $w$ ''downstairs'').

\end{problem}

Each solution of the problem defines a version of a {\it $W$-equivariant algebro-geometric distribution} ${\bf \Phi}_T$ on $T$ from \cite[Section 7.10]{BrKa} (the above condition (3) serves as a natural analogue of the requirement $3)$ from \cite[Section 7.10]{BrKa}). 

%In its turn, each such ${\bf \Phi}_T$ defines a {\it geometric lifting}  from $T$ to $G$ which, which, informally speaking, allows for replacing $T$ by %$G$ and the $W$-action by the adjoint $G$-action. 

Conjecture 7.11 from \cite{BrKa} asserts that for each algebraic $\ell$-dimensional representation  $\rho$ of the {\it Langlands dual group} $G^\vee$ there exists a $W$-equivariant algebro-geometric distribution ${\bf \Phi}_{\rho,T}$ with $X=\GG_m^\ell$  (where  $\GG_m$ stands for the multiplicative group), $f(c_1,\ldots,c_\ell)=\sum_{i=1}^\ell c_i$, and $\gamma=\gamma_\rho:X\to T$ is the homomorphism  of algebraic tori determined by $\rho$. The same conjecture claims that the existence of such ${\bf \Phi}_{\rho,T}$ implies a corollary  from  Local Langlands conjectures (see \cite[Section 1.1]{BrKa} for details).

Therefore, solving Problem \ref{prob:birational Weyl group actions} will help to dealing with the Local Langlands conjectures.

\begin{example} 
\label{ex:2 by 2}
Let $X=\GG_m^4$, $G=\{(A,A') \mid A, A' \in GL_2, \det A = \det A' \}$.

$T = \{(t_1, t_2; t'_1, t'_2) \mid t_1t_2 = t'_1t'_2 \}$.

$W = (\ZZ/2\ZZ)\times(\ZZ/2\ZZ)$.

$\gamma\begin{pmatrix}a & b \\ c & d\end{pmatrix} = (ab, cd; ac, bd) \in T,~ f \begin{pmatrix}a & b \\ c & d\end{pmatrix} = a+b+c+d \in \AA^1$.
$$
s_1\begin{pmatrix}a & b \\ c & d\end{pmatrix} = 
\begin{pmatrix} 
c\frac{a+d}{b+c} & d\frac{b+c}{a+d} \\
a\frac{b+c}{a+d} & b\frac{a+d}{b+c}
\end{pmatrix},~
s_2\begin{pmatrix}a & b \\ c & d\end{pmatrix} = \begin{pmatrix} 
b\frac{a+d}{b+c} & a\frac{b+c}{a+d} \\
 d\frac{b + c}{a+d} & c\frac{a+d}{b+c}
\end{pmatrix} \ .
$$
Clearly, $s_1s_2=s_2s_1$ and the function $f$ is $W$-invariant. It is also easy to see that $s_1(M)=M$  if and only if  $ab = cd$, i.e., fixed point of $s_1$ in $X=\GG_m^4$ are governed by the fixed points of $s_1$ in $T$. The same for $s_2$. One can show that the $W$-action on $X$ satisfying the requirements of Problem \ref{prob:birational Weyl group actions} is unique.

\end{example}

In order to solve Problem \ref{prob:birational Weyl group actions}, we introduced {\it geometric crystals} in \cite{bk1}. Let $I$ be the vertex set of the Dynkin diagram of $G$
and $\alpha^\vee_i : \GG_m \rightarrow T$, $\alpha_i: T \rightarrow \GG_m$ be respectively {\it simple coroots} and {\it simple roots} of $G$. The natural pairing between simple roots and simple coroots defines the Cartan matrix $a_{ij}=\left<\alpha_j, \alpha_i^\vee\right>$. In particular, 
if $G = GL_2$, 
$\alpha^\vee_1(c) = \begin{pmatrix}c & 0\\ 0 &c^{-1}\end{pmatrix}$, then 
$\alpha_1\begin{pmatrix}t_1& 0\\ 0 &t_2\end{pmatrix} = \displaystyle{\frac{t_1}{t_2}}$, $\left<\alpha_1, \alpha_1^\vee\right>=2$.

\begin{definition} A {\it decorated geometric crystal} is a $6$-tuple  ${\mathcal X}=(X,\gamma,f, \varphi_i,\varepsilon_i,e_i^\cdot|i\in I)$, where:

\noindent $\bullet$ $X$ is an irreducible algebraic variety.

\noindent $\bullet$ $\gamma$ is rational morphism $X\to T$. 

\noindent $\bullet$ $f,\varphi_i,\varepsilon_i:X\to \AA^1$ are rational functions. 

\noindent $\bullet$ each $e_i^\cdot:\GG_m\times X\to X$ is a unital rational
action of the multiplicative group
$\GG_m$ (to be denoted by $(c,x)\mapsto e_i^c(x)$) such that for each $i\in I$ one has:  

$$\gamma(e_i^c(x))=\alpha_i^\vee(c)\gamma(x),\varepsilon_i(x)=\alpha_i(\gamma(x))\varphi_i(x), \varepsilon_i(e_i^c(x))=c\varepsilon_i(x),\varphi_i(e_i^c(x))=c^{-1}\varphi_i(x)\ ,$$
\begin{equation}
\label{eq:f of e(x)}
f(e_i^c(x))=f(x)+
\frac{c-1}{\varphi_i(x)}+
\frac{c^{-1}-1}{\varepsilon_i(x)}
\end{equation}
for $x\in X$, $c\in \GG_m$; and for each $i\ne j$ one has the following geometric version of {\it Verma relations} (see \cite[Lemma 2.1]{bk1} and \cite[Proposition 39.3.7]{lu}):  

$\displaystyle{e_i^{c_1} e_j^{c_2}=e_j^{c_2} e_i^{c_1}}$
if $\left<\alpha_i,\alpha_j^\vee\right>=0$;

$\displaystyle{e_i^{c_1}e_j^{c_1c_2}e_i^{c_2}
=e_j^{c_2}e_i^{c_1c_2}e_j^{c_1}}$
if $\left<\alpha_j,\alpha_i^\vee\right>
=\left<\alpha_i,\alpha_j^\vee\right >=-1$;

$\displaystyle{e_i^{c_1}e_j^{c_1^2c_2}e_i^{c_1c_2}e_j^{c_2}
=e_j^{c_2}e_i^{c_1c_2}e_j^{c_1^2c_2}e_i^{c_1}}
$ if  $\left<\alpha_j,\alpha_i^\vee\right>=-2,
\left<\alpha_i, \alpha_j^\vee\right>=-1$;

$\displaystyle{
e_i^{c_1}
e_j^{c_1^3c_2}
e_i^{c_1^2c_2}
e_j^{c_1^3c_2^2}
e_i^{c_1c_2}
e_j^{c_2}
=e_j^{c_2}
e_i^{c_1c_2}
e_j^{c_1^3c_2^2}
e_i^{c_1^2c_2}
e_j^{c_1^3c_2}
e_i^{c_1}}
$
if $\left<\alpha_j,\alpha_i^\vee\right>=-3,
\left<\alpha_i,\alpha_j^\vee\right>=-1$.

\end{definition}

\begin{example} 
\label{ex:2 by 2 crystal}
In the notation of Example \ref{ex:2 by 2}, we have for $M=\begin{pmatrix} a & b \\ c & d\end{pmatrix}$: 
$$e^\tau_1(M) =\begin{pmatrix} 
\tau a\frac{a+d}{\tau a + d} & b\frac{\tau a + d}{a+d} \\
\tau^{-1} c\frac{\tau a + d}{a+d} & d\frac{a+d}{\tau a + d}
\end{pmatrix}\ , 
e^\tau_2(M) 
=\begin{pmatrix} 
\tau a\frac{a+d}{\tau a + d} & \tau^{-1}b\frac{\tau a + d}{a+d} \\
 c\frac{\tau a + d}{a+d} & d\frac{a+d}{\tau a + d}
\end{pmatrix} 
$$
and 
$$\varepsilon_1(M)=\frac{a + d}{cd},~\varphi_1(M)=\frac{a + d}{ab},~\varepsilon_2(M)=\frac{a + d}{bd},~\varphi_2(M)=\frac{a + d}{ac} \ .$$
Clearly, $e_1^\cdot$ commutes with $e_2^\cdot$ and both $f\circ e_1^\tau$ and $f\circ e_2^\tau$ satisfy \eqref{eq:f of e(x)}.

\end{example}

\begin{definition} For each decorated geometric crystal ${\mathcal X}=(X,\gamma,f, \varphi_i,\varepsilon_i,e_i^\cdot|i\in I)$ we define rational morphisms $s_i:X\to X$, $i\in I$ by 
$s_i(x):= e^{\frac{1}{\alpha_i(\gamma(x))}}_i(x)$.  

\end{definition}

It is easy to see that each $s_i$ is an involution. 
The following result shows that geometric crystals indeed solve Problem \ref{prob:birational Weyl group actions}.

\begin{proposition} For any decorated geometric crystal   ${\mathcal X}$ one has:

\noindent (a) The involutions $s_i$ satisfy the braid relations, i.e., define a rational action of $W$ on $X$.

\noindent (b) The function $f$ is  $s_i$-invariant for each $i\in I$. 

\end{proposition}

Part (a) coincides with \cite[Proposition 2.3]{bk1}, and we now prove part (b): 
$$f(s_i(x))=f(e^{\frac{1}{\alpha_i(\gamma(x))}}_i(x))=f(x)+
\frac{\frac{1}{\alpha_i(\gamma(x))}-1}{\varphi_i(x)}+
\frac{\alpha_i(\gamma(x))-1}{\varepsilon_i(x)}=f(x)$$
because $\varepsilon_i(x)=\alpha_i(\gamma(x))\varphi_i(x)$.

\begin{example} In the notation of Examples \ref{ex:2 by 2} and \ref{ex:2 by 2 crystal}, we have
$$e_1^{\frac{cd}{ab}} 
\begin{pmatrix} a & b \\ c & d\end{pmatrix} 
=s_1\begin{pmatrix}a & b \\ c & d\end{pmatrix},~  
e_2^{\frac{bd}{ac}} 
\begin{pmatrix} a & b \\ c & d\end{pmatrix}
=s_2\begin{pmatrix}a & b \\ c & d\end{pmatrix} \ .
$$

\end{example}

Now a new problem emerges: how to construct decorated geometric crystals. The answer comes from new geometric objects: linear unipotent bicrystals.

Let $U$ be a maximal unipotent subgroup of $G$ such that $TU=UT=B$ is a Borel subgroup, and let $\chi$ be a character of $U$, i.e., $\chi$ is a homomorphism $U\to \GG_a$ (where  $\GG_a$ stands for the additive group).

\begin{definition}

A {\it unipotent $\chi$-linear bicrystal} is a triple
$({\bf X},{\bf p},f)$ where:

\noindent $\bullet$   ${\bf X}$ is a  $U\times U$-{\it variety}, i.e., a pair
$(X,\alpha)$, where $X$ is an irreducible affine variety over $\QQ$ and
$\alpha:U\times X\times U\to X$ is a $U\times U$-action
on $X$, where the first $U$-action is {\it left} and second is {\it right},  such that each group $e\times U$ and $U\times e$ acts freely on $X$ (we will write the action as $(u,x,u')\mapsto uxu'$).

\noindent $\bullet$ ${\bf p}:X\to G$ is a $U\times U$-equivariant morphism, where the  action $U\times G\times U\to G$ is given by
$(u,g,u')\mapsto ugu'$.

\noindent $\bullet$ $f$ is a $\chi$-{\it linear} function on $X$, i.e., $
f(u\cdot x\cdot u')=\chi(u)+f(x)+\chi(u')$
for any $x\in X, u,u'\in  U$.

\end{definition}

The category of $\chi$-linear unipotent bicrystals is monoidal via the following {\it convolution product}:
$({\bf X},{\bf p},f)*({\bf Y},{\bf p}',f'):=({\bf X}*{\bf Y},{\bf p}'',f'')$, where:

\noindent $\bullet$ 
the $U\times U$-variety ${\bf X}*{\bf Y}$  is the quotient of $X\times Y$  by the following  left action of $U$ on $X \times Y$: $u\diamond(x,y)=(xu^{-1},u y)$.

\noindent $\bullet$  ${\bf p}'':X*Y\to G$ is defined by ${\bf p}''(x*y)={\bf p}(x){\bf p}'(y)$ 
for all $x\in X$, $y\in Y$.

\noindent $\bullet$ the function $f''$ on $Z=X*Y$ is defined by
$$f''(x*y)=f(x)+f'(y) $$
for all $x\in X$, $y\in Y$ clearly, both ${\bf p}''$ and $f$ are well-defined).

\begin{example} 
\label{ex:standard bicrystal}
Let $w_0\in W$ be the longest element of the Weyl group $W$ (i.e., the length of $w_0$ is $\dim U$). Take $X=Bw_0B$, the big Bruhat cell, ${\bf p}=id$ to be the natural inclusion $X\hookrightarrow G$, and $f_{G,\chi}:X\to \AA^1$ to be the function given by 
$$f_{G,\chi}(u\tilde w_0 u')=\chi(u)+\chi(u')$$
for any $u,u'\in U$ and any representative of $w_0$ of $W$ in the normalizer $Norm_G(T)$. Then the triple $(Bw_0B,id,f_{G,\chi})$ is a unipotent $\chi$-linear bicrystal.  

\end{example}

 Let $U_i$ be the one-parametric additive subgroup of $U$ corresponding to the simple root $\alpha_i$. And let  $\chi$ be a {\it regular} character, i.e., $\chi(U_i)\ne 0$ for all $i\in I$. For each $i\in I$ we choose a generator $x_i(a)$ of $U_i$ in such a way that $\chi(x_i(a))=a$ for $a\in \GG_a$. In particular, if $G=GL_n$, $U=U_n$, the group of upper uni-triangular matrices, $\chi(u)=\sum_{i=1}^{n-1} u_{i,i+1}$, and $x_i(a)=I+aE_{i,i+1}$.
 
\begin{example} 
\label{ex:explicit fG}
Assume that the group $G$ is simply-connected and $\chi$ is regular. Then in the notation of Example \ref{ex:standard bicrystal}, we have 
$$f_{G,\chi}(g)=\sum_{i\in I}\frac{\Delta_{w_0s_i\omega_i,\omega_i}(g)+\Delta_{w_0s_i\omega_i,\omega_i}(g)}{\Delta_{w_0\omega_i,\omega_i}(g)}$$
for each $g\in Bw_0B$,  where $\Delta_{\gamma,\delta}$ stands for a {\it generalized minor} defined in \cite{fz}. In particular, if $G=GL_n$, then 
$$f_{G,\chi}(g)=\sum_{i=1}^{n-1}\frac{\Delta_{\{n-i,n+2-i,\ldots,n\},\{1,\ldots,i\}}(g)+\Delta_{\{n+1-i,\ldots,n\},\{1,\ldots,i-1,i+1\}}(g)}{\Delta_{\{n+1-i,\ldots,n\},\{1,\ldots,i\}}(g)} \ ,$$
where $\Delta_{J,J'}(g)$ is the ordinary minor of an $n\times n$-matrix $g$ in the rows $J\subset \{1,2,\ldots,n\}$ and columns $J'\subset \{1,2,\ldots,n\}$. In particular, each denominator in the above formula is an $i\times i$ minor in the left lower corner which is never zero on $Bw_0B$. 

\end{example} 
 
We also fix $B^-$ to be the Borel subgroup opposite of $U$, i.e., $B^-\cap U=\{e\}$.  

Now we construct (decorated) geometric crystals out of $(U\times U,\chi)$-linear bicrystals:
\begin{equation}
\label{eq:functor F}
{\mathcal F}({\bf X}, {\bf p} ,f) := (X^-,\gamma,f^-,\varphi_i,\varepsilon_i,e_i^\cdot|i\in I)\ ,
\end{equation}
where:

\noindent $\bullet$ $X^-={\bf p}^{-1}(B^-)$. 

\noindent $\bullet$  $\gamma:X^-\to T$ is the composition of ${\bf p}:X^-\to B^-$ with the canonical projection $B^-\to B^-/U^-=T$.

\noindent $\bullet$ $f^-:X^-\to \AA^1$ is the restriction of the function $f$ to $X^-$.

\noindent $\bullet$ regular functions
$\varphi_i,\varepsilon_i:X^-\to \AA^1$, $i\in I$ are as follows. Let $pr_i$ be the   natural projection $B^-\to B^-\cap \phi_i(SL_2)$  (where $\phi_i$ is the $i$-th homomorphism $SL_2\to G$).
Using the fact that $x\in X^-$ if and only if ${\bf p}(x)\in B^-$, we set:
\begin{equation*}
\varphi_i(x):=\frac{b_{21}}{b_{11}},
%~ \alpha_i(x):=\frac{b_{11}}{b_{22}},
~ \varepsilon_i(x):=\frac{b_{21}}{b_{22}}=\varphi_i(x)\alpha_i(x)
\end{equation*}
for all $x\in X^-$, where $pr_i({\bf p}(x))=
\phi_i\begin{pmatrix}
b_{11} & 0 \\
b_{21} & b_{22}
\end{pmatrix}$.

\noindent $\bullet$ a rational morphism $e_i^\cdot:\GG_m \times X^- \to X$, $i\in I$ is given by ($x\in X^-$, $c\in \GG_m$):
\begin{equation}
\label{eq:simple multiplicative generator ei}
e_i^c(x)=x_i\left(\frac{c-1}{\varphi_i(x)}\right)\cdot x\cdot
x_i\left(\frac{c^{-1}-1}{\varepsilon_i(x)}\right) \ .
\end{equation}
if $\varphi_i\not = 0$ and $e_i^c(x)=x$ if $\varphi_i= 0$.

In particular, for $G=GL_2$, $X=Bw_0B$, one has $
e^c_1 \begin{pmatrix} b_{11} & 0 \\ b_{21} & b_{22} \end{pmatrix} 
= \begin{pmatrix} cb_{11} & 0 \\ b_{21} & c^{-1}b_{22} \end{pmatrix}$.

\begin{theorem} 
\label{th:from unipotent to geometric} ${\mathcal F}({\bf X}, {\bf p} ,f)$ is a decorated geometric crystal.
\end{theorem}

The ``non-decorated'' version of this result essentially coincides with \cite[Theorem 3.8]{bk1}. Let us demonstrate that  $f^-=f|_{X^-}$ satisfies \eqref{eq:f of e(x)}. Indeed, by \eqref{eq:simple multiplicative generator ei}, we have for each $x\in X^-$, $c\in \GG_m$, $i\in I$:
$$f^-(e_i^c(x))=f\left(x_i\left(\frac{c-1}{\varphi_i(x)}\right)\cdot x\cdot
x_i\left(\frac{c^{-1}-1}{\varepsilon_i(x)}\right)\right)$$
$$ =\chi\left(x_i\left(\frac{c-1}{\varphi_i(x)}\right)\right)
+f(x)+\chi\left(x_i\left(\frac{c^{-1}-1}{\varepsilon_i(x)}\right)\right)=f(x)+
\frac{c-1}{\varphi_i(x)}+
\frac{c^{-1}-1}{\varepsilon_i(x)} \ .$$

\section{Second lecture: Positive geometric crystals and crystal bases}

\medskip

We start with a combinatorial analogue of the geometric Problem \ref{prob:birational Weyl group actions} (see e.g., \cite{k93} for terminology on crystal bases).

\begin{problem}
Let $G^\vee$ be a complex reductive group, $X^\vee$ be a complex affine variety with $G^\vee$-action, $\mathcal{A}=\CC[X^\vee]$ be the coordinate algebra of $X^\vee$.
Construct a parametrization of the crystal basis $\BB$ for ${\mathcal A}=\CC[X^\vee]$  by a convex polyhedral cone in the lattice $\ZZ^\ell$ (where $\ell=\dim X^\vee$ in such a way that, under the parametrization, the functions $\tilde \varphi_i,\tilde \varepsilon_i$, and the crystal operators
$\tilde e_i^{\, n}:\ZZ^l \rightarrow \ZZ^l \sqcup \{\emptyset\}$ 
are given by piecewise linear formulas, i.e., using only $\min,\max,\pm$.

\end{problem}

All solutions $X^\vee$ we know so far (we will refer to them {\it good} varieties $X^\vee$) come from {\it positive} geometric crystals.

Some good $G^\vee$-varieties $X^\vee$ (positive geometric crystals with desirable properties exist):

\noindent $\bullet$ $X^\vee=G^\vee/U^\vee$ so that $\CC[X^\vee]=\oplus_\lambda V_{\lambda}$, where the sum is over all dominant weights of $G^\vee$.

\noindent $\bullet$  $G^\vee=\{(A,B)\in GL_m(\CC)\times GL_n(\CC),\det A=\det B\}$ and $X^\vee=\mbox{Mat}_{m\times n}(\CC)$ with the natural action
	$(A,B)(M)=A^{-1}MB$. More generally, $G^\vee$ is a Levi factor of a parabolic subgroup $P^\vee\subset \tilde G^\vee$, where $\tilde G^\vee$ is a larger complex reductive group and $X^\vee$ is the unipotent radical of $P^\vee$,  with the natural action of $G^\vee$ on $X^\vee$. 
	
\noindent $\bullet$ $G^\vee=\{(A,B,C)\in GL_2(\CC)\times GL_2(\CC)\times GL_2(\CC),\det A=\det B=\det C\}$ and $X^\vee=\mbox{Mat}_{2\times 2\times 2}(\CC)=\CC^2\otimes \CC^2\otimes \CC^2$ with the natural action of $G^\vee$.

A (conjecturally) bad $X^\vee$ is as follows: $G^\vee=GL_2(\CC)$, $X^\vee=S^3(\CC^2)$. It turns out that the geometric crystal with desirable properties does not exist.

Before describing positive geometric crystals let us define {\it positive varieties}. We first consider split algebraic tori $S= (\GG_m)^\ell$. 

\begin{definition} A positive morphism $f:(\GG_m)^k\to (\GG_m)^\ell$ is any rational morphism such that each coordinate function $f_i:(\GG_m)^k\to \GG_m$, $i=1,\ldots,\ell$ is a {\it positive} function, i.e., can be written as a ratio of two polynomials in $k$ variables with non-negative integer coefficients. In general, if $S$ and $S'$ are split algebraic tori (i.e., $S\cong (\GG_m)^k$ and $S'\cong (\GG_m)^\ell$), then a positive morphism $f:S \rightarrow S'$ is well-defined.

%$\forall \lambda^\vee,\mu, \GG_m \stackrel{\lambda^\vee}{\rightarrow} S 
%\stackrel{f}{\rightarrow} S' \stackrel{\mu}{\rightarrow} \GG_m$ 
%is a well defined {\bf positive} function in one variable.

\end{definition}

For instance, e.g. if $S=\GG_m$, $f(x)=\frac{x^3+1}{x+1}=x^2-x+1$ is positive. Note also that the morphism  $f:(\GG_m)^2 \rightarrow (\GG_m)^2$ given by $f(x,y)=(x,x+y)$ is positive but its inverse $f^{-1}(x,y)=(x,y-x)$ is not positive.

It is easy to see that composition of positive morphisms are positive. Consider a category $\mathcal{T}_+$ whose objects are split algebraic tori 
and arrows are positive morphisms.
%Let $S$ be a split algebraic torus i.e., $S\cong \GG_m^l$. Let 
%$X_*(S)=\mbox{Hom}(\GG_m,S)$ be the lattice of co-character,
%$X^*(S)=\mbox{Hom}(S,\GG_m)$ be the lattice of character.

Now we construct a ``tropicalization'' functor $Trop:\mathcal{T}_+ \to {\bf Sets}$ as follows. First, for a non-zero Laurent polynomial $f(x)=\sum_{i=n}^N a_i x^i$ with $a_n\ne 0$, $N\ge n$, we set $\deg f=n$. And for any rational function $f:\GG_m\to \GG_m$ written as the ratio of two Laurent polynomials $f=\frac{f_1}{f_2}$, we set $\deg f=\deg f_1-\deg f_2$. Then we set $Trop(S)=X_\star(S)=Hom(\GG_m,S)$, the lattice of co-characters of a split algebraic torus $S$; and for each positive morphism $f:S\to S'$ we set $Trop(f)$ to be a (piecewise-linear) map $X_\star(S)\to X_\star(S')$ determined by:
$$\langle\mu,Trop(f)(\lambda)\rangle= \deg f_{\lambda,\mu}\ ,$$
for any co-character $\lambda\in X_\star(S)$ and any character $\mu\in X^\star(S)=Hom(S,\GG_m)$, where $f_{\lambda,\mu}:\GG_m \stackrel{\lambda}{\to}
S\stackrel{f}{\rightarrow}S'\stackrel{\mu}{\rightarrow}\GG_m$, and $\langle\bullet ,\bullet \rangle:X^\star(S)\times X_\star(S)\to \ZZ$ is the canonical pairing of characters and co-characters.

Therefore, for each rational function $f(x_1,\ldots,x_k)=(\sum_{a\in \ZZ^k} c_a x^a)/(\sum_{a\in \ZZ^k} d_a x^a)$, where we abbreviated  $x^a=x_1^{a_1}\cdots x_k^{a_k}$,  the tropicalization $Trop(f):\ZZ^k\to \ZZ$ is a piecewise-linear function given by: 
$$Trop(f)(\tilde x)=\min_{a\in \ZZ^k:c_a\ne 0}\left (\sum_{i=1}^k a_i\tilde x_i\right )-\min_{a\in \ZZ^k:d_a\ne 0}\left (\sum_{i=1}^k a_i\tilde x_i\right )\ ,$$ 
and for each rational morphism $f:(\GG_m)^k\to (\GG_m)^\ell$ given by $f(x)=(f_1(x),\ldots,f_\ell(x))$, the tropicalization $Trop(f):\ZZ^k\to \ZZ^\ell$ is a piecewise-linear map given by 
$$Trop(f)(\tilde x)=(Trop(f_1)(\tilde x),\ldots,Trop(f_\ell)(\tilde x)) \ .$$

\begin{example} If $f(x)=\frac{x^3+1}{x+1}=x^2-x+1$, then $Trop(f):\ZZ\to \ZZ$ is given by
$$Trop(f)(\tilde x)=\min(3\tilde x,0)-\min(\tilde x,0)=\min(2\tilde x,\tilde x,0)=\min(2\tilde x,0) \ .$$
And if $f:(\GG_2)^2\to (\GG_2)^2$ is given by $f(x,y)=(x,x+y)$, then 
$$Trop(f)(\tilde x,\tilde y)=(\tilde x,\min(\tilde x,\tilde y)) \ .$$
\end{example}

\begin{theorem} 
\label{th:Trop}
\cite[Corollary 2.10]{bk1} $Trop:\mathcal{T}_+\to {\bf Sets}$ is, indeed, a functor.
\end{theorem}

Note that the positivity is important for the functoriality - in the above example, the morphism $f:(\GG_2)^2\to (\GG_2)^2$ is positive and  invertible, but its inverse, given by $f^{-1}(x,y)=(x,y-x)$, is not positive. It is easy to see that $Trop(f^{-1})\circ Trop(f)\ne Id$.

\begin{definition} A {\it positive variety} is a pair $(X,\theta)$, where $X$ is an irreducible variety (defined over $\QQ$) and $\theta:S\widetilde \to X$ is a birational isomorphism of a split algebraic torus $S$ and $X$. A morphism of positive varieties $(X,\theta)\to (Y,\theta')$ or {\it $(\theta,\theta')$-positive morphism $X\to Y$} (where $\theta:S\widetilde \to X$ and $\theta':S'\widetilde \to X$) is a rational morphism $f:X\to Y$ such that ${\theta'}^{-1}\circ f\circ \theta$ is a well-defined positive morphism $S\to S'$.

\end{definition}

\begin{proposition} 
\label{pr:monoidal positive variety}
Positive varieties and their morphisms form a category. This category is monoidal with respect to the product $(X,\theta)\times (Y,\theta'):=(X\times Y,\theta\times \theta')$. 
\end{proposition}

\begin{proof} Indeed, if $\theta:S\widetilde \to X$ and $\theta':S\widetilde \to Y$ are birational isomorphisms, then so is $\theta\times \theta':S\times S'\to X\times Y$. That is, the product is well-defined. Its associativity follows, and the unit object is the pair $(S_0,id)$, where $S_0=\{e\}$ is the $0$-dimensional torus.
\end{proof}

By definition,  for each split algebraic torus $S$ the pair $(S,id)$ is a natural positive variety and $(S,id)\to (S',id)$ is a natural morphism of positive varieties.

We  say that for a positive variety $(X,\theta)$  a non-zero function $f:X\to \AA^1$ is $\theta$-positive if  $f:(X,\theta)\to (\AA^1,id)$ is a morphism of positive varieties.
 
\begin{definition}  A {\it positive decorated geometric crystal} is a pair $({\mathcal X},\theta)$, where ${\mathcal X}=(X,\gamma,f, \varphi_i,\varepsilon_i,e_i^\cdot|i\in I)$ is a geometric crystal, $(X,\theta)$ is a positive variety, and:

\noindent $\bullet$ $\gamma:(X,\theta) \rightarrow (T,id)$ is a morphism of positive varieties.

\noindent $\bullet$ The function $f:X \rightarrow \AA^1$ is $\theta$-positive.

\noindent $\bullet$ The action $e_i^\cdot$ is a morphism of positive varieties $(\GG_m,id)\times (X,\theta)\to  (X,\theta)$ is a morphism of positive varieties.

\end{definition}

Note that if $({\mathcal X},\theta)$ is positive, then the functions $\varepsilon_i,\varphi_i$ are also positive.

For each positive variety $(X,\theta)$ (where $\theta:S\widetilde \to X$) we denote $Trop(X,\theta):=Trop(S)=X_\star(S)$. And for each morphism $f:(X,\theta)\to (Y,\theta')$ of positive varieties we denote $Trop(f):=Trop({\theta'}^{-1}\circ f\circ \theta):(X,\theta)\to Trop(Y,\theta')$.

\begin{lemma} The association $(X,\theta)\mapsto Trop(X,\theta)$ is a monoidal functor from the category of positive varieties to the category of sets. 

\end{lemma}   

\begin{proof} Indeed, the functoriality follows from Theorem \ref{th:Trop}. And the fact that the functor respects products follows from Proposition \ref{pr:monoidal positive variety}. \end{proof} 

Now apply the functor $Trop$ to  positive geometric crystals. For any positive decorated geometric crystal $({\mathcal X},\theta)$  denote $Trop({\mathcal X},\theta)=(\widetilde X,\tilde \gamma,\tilde f,\tilde \varphi_i,\tilde \varepsilon_i,\tilde e_i^\cdot| i\in I)$, where:

\noindent $\bullet$  $\tilde X=Trop(X,\theta)$ is the ''tropical'' variety.

\noindent $\bullet$ $\tilde \gamma=Trop(\gamma):\tilde X\to X_*(T)=\Lambda^\vee$ is the co-weight grading.

\noindent $\bullet$ $\tilde e_i^\cdot=Trop(e_i^\cdot):\ZZ\times \tilde X\to \tilde X$ is a free $\ZZ$-action

\noindent $\bullet$ $\tilde f=Trop(f):\tilde X\to \ZZ$. 
 
\noindent $\bullet$ $\tilde{\varepsilon}_i,\tilde{\varphi}_i$ are functions $\tilde X\to \ZZ$  given by $\tilde{\varepsilon}_i=-Trop(\varepsilon_i)$, $\tilde{\varphi}_i=-Trop(\varphi_i)$.

Clearly, each $e_i^\cdot$ defines two mutually inverse bijections $\tilde e_i:\tilde X\to \tilde X$ and $\tilde f_i:\tilde X\to \tilde X$ via $\tilde e_i=\tilde e_i^{\,1}$ and $\tilde f_i=\tilde e_i^{\,-1}$.

\begin{theorem} $Trop({\mathcal X},\theta)$ is a torsion-free  Kashiwara crystal such that:
\begin{equation}
\label{eq:f of tilde e decorated}
\tilde f(\tilde e_i^{\,n}(\tilde b))=\min(\tilde f_0(\tilde b),n+\tilde \varphi_i(\tilde b),-n+\tilde \varepsilon_i(\tilde b))
\end{equation}
for $\tilde b\in \widetilde{X}$, $n\in \ZZ$, $i\in I$, where $\tilde f_0:\tilde X\to \ZZ$ is a function.

\end{theorem}

The ``non-decorated'' version of this result coincides with \cite[Theorem 2.11]{bk1}. To prove \eqref{eq:f of tilde e decorated}, let us rewrite \eqref{eq:f of e(x)} as follows:
\label{cl:f of e}
\begin{equation}
\label{eq:f of e decorated}
f(e_i^c(x))=f_0(x)+
\frac{c}{\varphi_i(x)}+
\frac{c^{-1}}{\varepsilon_i(x)}
\end{equation}
for $x\in X$, $c\in \GG_m$, $i\in I$, where $f_0(x)=f(x)-\frac{1}{\varphi_i(x)}-\frac{1}{\varepsilon_i(x)}$. Denote  $\tilde f_0:=Trop(f_0)$ (if $f_0= 0$, then $\tilde f_0=+\infty$). Note that the  function $\GG_m\times X\to \AA^1$ given by $(c,x)\mapsto f(e_i^c(x))$ is positive. It follows from  Theorem \ref{th:Trop} that the tropicalization of this positive function is the function  $\ZZ\times \widetilde{X}\to \ZZ$ given by $(n,\tilde b)\mapsto \tilde f(\tilde e_i^{\,n}(\tilde b))$. Therefore, applying tropicalization functor to the identity \eqref{eq:f of e decorated}, we obtain
%$$\tilde f(\tilde e_i^{\,n}(b))=\min(\tilde f_0(b),\tilde f_1(n,b)) \ ,$$
%and taking into account that $\tilde f_1(n,b)=\min(n+\tilde \varphi_i(b),-n+\tilde \varepsilon_i(b))$ we obtain:
\eqref{eq:f of tilde e decorated}.
\endproof

Note that the Kashiwara crystal $Trop({\mathcal X},\theta)$ is associated to the Langlands dual group $G^\vee$ rather than $G$ because it is graded by the co-weights of $G$, i.e., the weights of $G^\vee$. 

Let $\tilde B:=\{\tilde b\in Trop(X,\theta) \mid \tilde{f}(\tilde b)\ge 0\}$ and denote by $\BB({\mathcal X},\theta)$ the restriction of the free Kashiwara crystal $Trop({\mathcal X},\theta)$ to $\tilde B$.

Note that if the function $f\circ \theta:S\to \AA^1$ is regular, then $\tilde B$ is a convex
polyhedral cone in $\tilde X$.

\begin{proposition} The Kashiwara crystal  $\BB({\mathcal X},\theta)$ is normal.
\end{proposition}

\begin{proof} Recall from \cite{k93} that the normality condition is that:
\begin{equation}
\label{eq:normality}
\varepsilon_i(\tilde b)=\max\{n\ge 0:\tilde e_i^{\,n}(\tilde b)\ne \emptyset\},~\varphi_i(\tilde b)=\max\{n\ge 0:\tilde e_i^{\,-n}(\tilde b)\ne \emptyset\} \ .
\end{equation}
According to \eqref{eq:f of tilde e decorated},   $\tilde e_i^{\,n}(\tilde b)\ne \emptyset$ for $\tilde b\in \tilde B$ if and only if $-\tilde \varphi_i(\tilde b)\le n\le \varepsilon_i(\tilde b)$. This is exactly the normality condition \eqref{eq:normality}.
\end{proof}

Next, we construct our main example of a positive geometric crystal $({\mathcal X},\theta)$ such that $\BB({\mathcal X},\theta)$ is a crystal basis for an integrable $G^\vee$-module. 

For each sequence $\ii=(i_1,\ldots,i_\ell)\in I^\ell$  we define a  morphism 
$\theta_\ii^-:T\times (\GG_m)^\ell\widetilde\to B^-$ by the formula:
\begin{equation}
\label{eq:thetaii}
\theta_\ii(t,c_1,\ldots,c_\ell):=t\cdot x_{-i_1}(c_1)\cdot x_{-i_2}(c_2)\cdots x_{-i_\ell}(c_\ell)
\end{equation}
for any $c_1,\ldots,c_\ell\in \GG_m$,  where
$x_{-i}:\GG_m\to B^-$ is given by the formula
\begin{equation}
\label{eq:thetai}
x_{-i}(c):= \begin{pmatrix}
c^{-1}& 0 \\
1 & c
\end{pmatrix}_i
\end{equation}
where $g\mapsto g_i$ is the homomorphism $SL_2\to G$ such that 
$\begin{pmatrix}
c& 0 \\
0 & c^{-1}
\end{pmatrix}_i=\alpha_i^\vee(c)$ and $\begin{pmatrix}
1& a \\
0 & 1
\end{pmatrix}_i=x_i(a)$.

Now let ${\mathcal X}_G:={\mathcal F}(Bw_0B,id,f_{G,\chi})$ be the decorated geometric crystal associated to the $(U\times U,\chi)$-linear bicrystal from Example \ref{ex:standard bicrystal}, where $\chi$ is a regular character of $U$. Note that the underlying variety of ${\mathcal X}_G$ is $(Bw_0B)^-=Bw_0B\cap B^-$.

\begin{proposition}  
\label{pr:reduced word}
Let $\ii=(i_1,\ldots,i_\ell)$ be a reduced decomposition of the longest element $w_0\in W$ (i.e., $\ell=\dim U$ and $s_{i_1}\cdots s_{i_\ell}=w_0$). Then: 

\noindent (a) The morphism $\theta_\ii$ is an open embedding (and hence a birational isomorphism) $T\times (\GG_m)^\ell\hookrightarrow B^-\cap Bw_0B$.

\noindent (b) The pair $({\mathcal X}_G,\theta_\ii)$ is a positive decorated geometric crystal. 
\end{proposition} 

Part (a) of Proposition \ref{pr:reduced word} follows from \cite[Theorems 1.2, 1.3]{fz}, and part (b) follows from the results of \cite[Section 3.2]{bk2}.

By tropicalizing this positive geometric crystal, we obtain our main result.
 
\begin{theorem} \cite[Main Theorem 6.15]{bk2}. For any reduced decomposition $\ii$ of $w_0$ the normal Kashiwara crystal $\BB({\mathcal X}_G,\theta_\ii)$ is isomorphic to the disjoint union of all irreducible $G^\vee$-crystal bases ${\mathcal B}_\lambda$. 
\end{theorem}

The proof of this result is rather non-trivial. It is based on the notion of {\it strongly positive} $\chi$-linear unipotent bicrystals introduced in \cite[Section 3.2]{bk2} and Joseph's characterization of the irreducible crystal bases via {\it closed families} (\cite[Section 6.4.21]{Joseph}). 

\begin{example} Let $G=GL_3$, so that $T=\{t=diag(t_1,t_2,t_3)\}\subset GL_3$. We fix the reduced decomposition $\ii=(1,2,1)$ of $w_0\in W=S_3$ so that: 
$$\theta_\ii(t;c_1,c_2,c_3)=\begin{pmatrix}
t_1&0 &  0\\
0& t_2&0\\
0&0&t_3\end{pmatrix} 
\begin{pmatrix}
c_1^{-1}&0 &  0\\
1& c_1&0\\
0&0&1\end{pmatrix}
\begin{pmatrix}
1&0 &  0\\
0& c_2^{-1}&0\\
0&1&c_2\end{pmatrix}
\begin{pmatrix}
c_3^{-1}&0 &  0\\
1& c_3&0\\
0&0&1\end{pmatrix}=$$
$$=\begin{pmatrix}
t_1\frac{1}{c_1c_3}&0 & 0\\
t_2(\frac{c_1}{c_2}+\frac{1}{c_3})& t_2\frac{c_1c_3}{c_2}&0\\
t_3&t_3c_3&t_3c_2
\end{pmatrix} \ .$$ 

Therefore, according to Example \ref{ex:explicit fG}, the restriction of $f_{G,\chi}$ to $B^-\cap Bw_0B$ is given by (in the new coordinates $(t;c_1,c_2,c_3)$):
$$f_{G,\chi}(t;c_1,c_2,c_3)=c_1+\frac{c_2}{c_3}+c_3+\frac{t_2}{t_3}\cdot \left (\frac{c_1}{c_2}+\frac{1}{c_3}\right )+\frac{t_1}{t_2}\cdot \frac{1}{c_1} \ .$$
And the rest of the decorated geometric crystal structure ${\mathcal X}_G$ on $B^-\cap Bw_0B$ is given by the morphism $\gamma$, the actions $e_i^\cdot$, and the functions $\varphi_i,\varepsilon_i$, $i=1,2$:
$$\gamma(t;c_1,c_2,c_3)=\left (t_1\frac{1}{c_1c_3},t_2\frac{c_1c_3}{c_2},t_3c_2\right )\ ,$$ 
$$e_1^d(t;c_1,c_2,c_3)=\left(t;c_1\frac{c_2+c_1c_3}{d\cdot c_2+c_1c_3},c_2,c_3\frac{c_2+d^{-1}\cdot c_1c_3}{c_2+c_1c_3}\right) \ ,$$
$$e_2^d(t;c_1,c_2,c_3)=(t;c_1,d^{-1}\cdot c_2,c_3) \ .$$
$$\varphi_1(t;c_1,c_2,c_3)=\frac{t_2}{t_1}\cdot \left(\frac{c_1^2c_3}{c_2}+c_1\right), \varphi_2(t;c_1,c_2,c_3)=\frac{t_3}{t_2}\cdot \frac{c_2}{c_1}\ ,$$
$$\varepsilon_1(t;c_1,c_2,c_3)= \frac{1}{c_3}+\frac{c_2}{c_1c_3^2},~ \varepsilon_2(t;c_1,c_2,c_3)=\frac{c_3}{c_2} \ .$$

The tropicalization of the above structures consists of:

\noindent $\bullet$ The set  $\tilde X=\Lambda^\vee\times \ZZ^3$ where $\lambda\in \Lambda^\vee= \ZZ^3$. 

\noindent $\bullet$ The functions $\tilde f_{G,\chi},\tilde \varphi_i,\tilde \varepsilon_i:\Lambda^\vee \times \ZZ^3\to \ZZ$, $i=1,2$:
$$\tilde f_{G,\chi}(\lambda;{\bf m})=\min(m_1,m_2-m_3,m_3,\lambda_2-\lambda_3-\max(m_3,m_2-m_1),\lambda_1-\lambda_2-m_1)\ ,$$
$$\tilde \varphi_1(\lambda;{\bf m})=\lambda_1-\lambda_2-\min(m_1,2m_1+m_3-m_2),~\tilde \varphi_2(\lambda;{\bf m})=\lambda_2-\lambda_3+m_1-m_2\ ,$$
$$\tilde \varepsilon_1(\lambda;{\bf m}))= \max(m_3,m_1+2m_3-m_2),~\tilde \varepsilon_2(\lambda;{\bf m})=m_2-m_3 $$
for $(\lambda;{\bf m})=((\lambda_1,\lambda_2,\lambda_3);(m_1,m_2,m_3))\in \tilde X$.

\noindent $\bullet$ The set $\tilde B=\{(\lambda;{\bf m})\in \tilde X:\tilde f_{G,\chi}(\lambda;{\bf m})\ge 0\}$, i.e., $\tilde B$ consists of all $(\lambda;{\bf m})\in \tilde X$ such that $m_1\ge 0,~m_2\ge m_3\ge 0,~\lambda_1-\lambda_2\ge m_1,~\lambda_2-\lambda_3\ge m_3, \lambda_2-\lambda_3\ge m_2-m_1$.
That is, each point of $(\lambda;{\bf m})\in \tilde B$ is a Gelfand-Tsetlin pattern:
$$
\begin{pmatrix}
\lambda_1 &               & \lambda_2    &              & \lambda_3\\
          & \lambda_2+m_1 &              &\lambda_2+m_3 & \\
          &               &\lambda_3+m_2 &              &
\end{pmatrix}
$$

\noindent $\bullet$ The bijection $\tilde e_i^{\,n}:\tilde B\to \tilde B$, $i=1,2$:
$$\tilde e_1^{\,n}(\lambda;{\bf m})=\left(\lambda;m_1+\max(\delta-n,0)-\max(\delta,0),m_2,m_3+\max(\delta,0)-\max(\delta,n)\right)\ ,$$
where $\delta=m_1+m_3-m_2$,
$$\tilde e_2^{\,n}(\lambda;{\bf m})=\left(\lambda;m_1,m_2-n,m_3\right) \ .$$
Therefore, one has the decomposition into the connected components. 
$${\mathcal B}({\mathcal X}_G,\theta_\ii)\cong \bigsqcup_{\lambda=(\lambda_1\ge \lambda_2\ge \lambda_3)} {\mathcal B}_\lambda \ .$$

\end{example}

\end{document}